\documentclass[12pt]{article}
\usepackage[utf8]{inputenc}

\usepackage[T1]{fontenc}
\usepackage{lmodern}
\usepackage{amsmath,amssymb,amsthm}
\usepackage{blindtext}
\usepackage{graphicx}

\usepackage[all]{xy}

\usepackage{enumitem}
\usepackage{xcolor}
\definecolor{darkcandyapplered}{rgb}{0.64, 0.0, 0.0}
\definecolor{joliblanc}{RGB}{252, 251, 249}
\definecolor{bluegray}{rgb}{0.4, 0.6, 0.8}
\definecolor{ceruleanblue}{rgb}{0.16, 0.32, 0.75}
\definecolor{greydark}{RGB}{75, 71, 59}
\definecolor{oldlavender}{RGB}{75, 71, 59}

\colorlet{shadecolor}{white}

\usepackage{hyperref}
\hypersetup{colorlinks=true,linkcolor=darkcandyapplered,citecolor=oldlavender}

\usepackage{caption}

\def\defn#1{{\itshape #1}}
\def\cit#1{{ \scshape #1}}

\usepackage{thmtools}

\declaretheoremstyle[%
  spaceabove=-6pt,%
  spacebelow=6pt,%
  headfont=\normalfont\itshape,%
  postheadspace=1em,%
  qed=\qedsymbol,%
  headpunct={}
]{mystyle} 
\declaretheorem[name={Proof ---},style=mystyle,unnumbered,
]{Proof}

\renewcommand{\qedsymbol}{\color{oldlavender}{$\blacksquare$}}

\usepackage{lipsum}
\usepackage{amsthm}
\usepackage{etoolbox}
\usepackage{changepage}
\usepackage{framed}

\newtheoremstyle{mytheoremstyle} 
    {\topsep}                    
    {\topsep}                    
    {\itshape}                   
    {}                           
    {\scshape}                   
    {.}                          
    {.5em}                       
    {}  

\theoremstyle{mytheoremstyle}

\newenvironment{frshaded*}{%
    \MakeFramed {\advance\hsize-\width \FrameRestore}}%
    {\endMakeFramed} 
    
\newtheorem{theo}{Theorem}[section]
\newtheorem{prop}[theo]{Proposition}
\newtheorem{corollaire}[theo]{Corollary}
\newtheorem{lemme}[theo]{Lemma}
\newtheorem{conj}[theo]{Conjecture}

\newtheoremstyle{mydefinitionstyle} 
    {\topsep}                    
    {\topsep}                    
    {}                   
    {}                           
    {\itshape}                   
    {.}                          
    {.5em}                       
    {}  

\theoremstyle{mydefinitionstyle}

\newtheorem{definition}{Definition}[section]
\newtheorem*{remark}{Note}
\newtheorem{example}[definition]{Example}

\newenvironment{theorem}
      {\setlength\FrameRule{4 pt} \begin{frshaded*}\begin{theo}}
      {\end{theo}\end{frshaded*}}

\newenvironment{lemma}
      {\setlength\FrameRule{4 pt}\begin{frshaded*}\begin{lemme}}
      {\end{lemme}\end{frshaded*}}

\newenvironment{cor}
      {\setlength\FrameRule{4 pt}\begin{frshaded*}\begin{corollaire}}
      {\end{corollaire}\end{frshaded*}}

\newenvironment{proposition}
      {\setlength\FrameRule{4 pt}\begin{frshaded*}\begin{prop}}
      {\end{prop}\end{frshaded*}}

\usepackage{geometry}

\usepackage{geometry}
\usepackage{titlesec}
\usepackage{secdot}
\titleformat*{\section}{\centering\bfseries}
\titleformat*{\subsection}{\itshape}

\usepackage{enumitem}

\usepackage{bbm}

\newcommand{\T}{\mbox{T}}

\newcommand{\R}{\mathbb{R}}
\newcommand{\Z}{\mathbb{Z}}
\newcommand{\N}{\mathbb{N}}
\newcommand{\C}{\mathbb{C}}

\newcommand{\J}{\Phi}

\newcommand\restr[2]{{
  \left.\kern-\nulldelimiterspace 
  #1 
  \vphantom{\big|} 
  \right|_{#2} 
  }}

\newcommand{\cat}{\mbox{Cat}}
  
\usepackage{authblk}
\author{\small\textsc{Marine Fontaine} and \textsc{James Montaldi}}

\title{A localization formula for equivariant Lyusternik-Schnirelmann category}
\date{}
\begin{document}
\maketitle

\begin{abstract}
	 The LS-category of a topological space is a numerical homotopy invariant, introduced originally in a course on the global calculus of variations by Lyusternik and Schnirelmann, to estimate the number of critical points of a smooth function. When the topological space is a smooth manifold equipped with a proper action of a Lie group, we give a localization formula to calculate the equivariant analogue of this category in terms of the minimal orbit-type strata. The formula holds provided that the manifold admits a specific cover. We show that such a cover exists on every symplectic toric manifold. The known result stating that the LS-category of a symplectic toric manifold is equal to the number of fixed points of the torus action follows from our localization formula.
\end{abstract}

\tableofcontents
\section{Introduction}
The \defn{Lyusternik-Schnirelmann category} or \defn{LS-category} of a topological space $M$ is the numerical homotopy invariant $\cat(M)$ defined to be the least number of open subsets $U\subset M$, whose inclusion is nullhomotopic, that are required to cover $M$. Although it is now the subject of a full theory in connection with algebraic topology, it was originally introduced by Lyusternik and Schnirelmann in a course on the global calculus of variations, when $M$ is a smooth closed manifold without boundary~\cite{LS}. In this case they show that any smooth real-valued function $f$ defined on $M$ has at least $\cat(M)$ critical points. The difference with Morse theory is that $f$ is allowed to have degenerate critical points. Rewiews on the Lyusternik-Schnirelmann theory are for instance~\cite{James,MR1990857,bartsch}. 

If $M$ is a topological space, acted on continuously and properly by a topological group $G$, then the LS-category has an equivariant analogue $\cat_G(M)$. It has first been introduced by\cit{Fadell}~\cite{fadell} and\cit{Marzantowicz}~\cite{marzantowicz} for compact groups, and by\cit{Colman}~\cite{colman} for finite groups. A substantial part of the theory has been extended for proper group actions by\cit{Ayala, Lasheras and Quintero}~\cite{ayala}.

A \defn{$G$-categorical} open subset of $M$ is a $G$-invariant open subset $U\subset M$ admitting a $G$-deformation retract onto a $G$-orbit (cf. Definition \ref{def: G-categorical open subset}). The numerical homotopy invariant $\cat_G(M)$ is the least number (possibly infinite) of $G$-categorical open subsets that are required to cover $M$. If $M$ is a smooth manifold and $G$ is a Lie group, a class of $G$-categorical open subsets consists of \defn{$G$-tubular} open subsets, which are essentially tubular neighbourhoods of group orbits (cf. Definition \ref{def:G tubular open subset}). This fact is a direct application of the Tube Theorem \ref{Tube theorem}.

This invariant is in general difficult to compute and we are usually only able to know an estimation of it, in term of the cup length of $M$. We obtain a new formula to reduce the calculation of $\cat_G\left(M\right)$ to the calculation of the equivariant LS-category of the minimal orbit-type strata of $M$. In general, any topological space $M$ can be written as a disjoint union of smaller subsets $M_{\beta}$, called \defn{strata}, indexed on some strictly partially ordered set $(\mathcal{B},\prec)$. Those strata are required to fit in a specific way and form themselves a strictly partially ordered set (cf. Section \ref{isotropytype}). A stratum is minimal if it is minimal with respect to the strict partial order defined on them. If $M$ is a proper $G$-manifold, the strata $M_{\beta}$ are generally the connected components of the orbit-type submanifolds. We use a modified definition of orbit-type stratum (cf. Definition \ref{def: orbit type stratum}). We say that an \defn{orbit-type stratum} is a $G$-orbit of a connected component of the subset of $M$ of all the points having the same stabilizer.

On a large class of proper $G$-manifolds $M$, including symplectic toric manifolds, we observe that $M$ can be entirely covered by a subcover of its minimal orbit-type strata, made of $G$-tubular open subsets. Besides this cover is the smallest cover, made of $G$-categorical open subsets, that we can take. Such covers are called \defn{minimal $G$-tubular covers} and are discussed in Section \ref{s:G-categorical tubular covers}. However those covers do not exist in general. A non-example is given, when $M$ is a non-Hamiltonian compact $S^1$-manifold (cf. Example \ref{non-example}). By using the natural stratification of the moment polytope, we show in Section \ref{s: tubular toric} that every symplectic toric manifold admits a minimal $G$-tubular cover, where $G$ in this case is a torus having half the dimension of $M$ and acting effectively on it (cf. Theorem \ref{thm toric}). When $M$ admits a minimal $G$-tubular cover, we show that the calculation of $\cat_G(M)$ is intrinsically reduced to those of the minimal orbit-type strata of $M$. Explicitly, we obtain the localization formula
\begin{equation*}
	\cat_G(M)=\sum\cat_{G}\left(M_{\beta}\right)
\end{equation*}
where the summation is taken over the minimal orbit-type strata $M_{\beta}$ (cf. Section \ref{s: loc formula}, Theorem \ref{localization thm first} and Corollary \ref{localization thm}). The result of\cit{Bayeh and Sarkar} (cf. \cite{MR3422738} Theorem $5.1$), which states that the equivariant Lyusternik-Schnirelmann category of a quasitoric manifold is precisely the number of fixed points of the torus action, is a consequence of Theorem \ref{thm toric} and of our localization formula.
\paragraph{Acknowledgements.} We would like to thank Yael Karshon and Eckhard Meinrenken for useful discussions and suggestions.
\vspace{1cm}

After completing this work, which forms part of the thesis \cite{moi}, we discoverd the paper of \cit{Hurder and T\"oben} \cite{MR2492297} which contains one of our main theorems (Theorem \ref{localization thm first}) and uses a similar approach.
\section{Terminologies}
We work with smooth manifolds and, except stated otherwise, the term \defn{submanifold} refers to an embedded submanifold. A smooth \defn{action} of a Lie group $G$ on $M$ is a group homomorphism $G\to \mbox{Diff}(M)$. We denote the action map by $(g,m)\in G\times M\mapsto g\cdot m\in M$. A \defn{$G$-manifold} is a pair $(M,G)$ where $M$ is a smooth manifold acted on by a Lie group $G$. If the action is proper then we refer to $(M,G)$ as a \defn{proper $G$-manifold}. A smooth map $f:M\to N$ between two $G$-manifolds is \defn{$G$-equivariant} if $f(g\cdot m)=g\cdot f(m)$ for all $g\in G$ and $m\in M$. The \defn{stabilizer} of $m\in M$ is the subgroup $G_m=\lbrace g\in G\mid g\cdot m=m\rbrace$. We say that the action of $G$ on $M$ is \defn{free} if all the stabilizers $G_m$ are equal to the trivial group $\mathbbm{1}$. The \defn{group orbit} or \defn{$G$-orbit} of a point $m\in M$ is the set $G\cdot m=\lbrace g\cdot m\mid g \in G\rbrace$. We use the notation $\mathfrak{g}\cdot m$ to mean the tangent space to $G\cdot m$ at $m$. The space $\mathfrak{g}$ stands for the Lie algebra of $G$ with Lie bracket $[\cdot,\cdot ]$, obtained by identifying $\mathfrak{g}$ with the left invariant vector fields on $G$.
\subsection*{Local models of G-orbits}
Let $(M,G)$ be a proper $G$-manifold. In this case, the group orbits are embedded submanifolds and all the stabilizers $G_m$ are compact. Let $N\subset T_mM$ be a $G_m$-invariant complement to $\mathfrak{g}\cdot m$ in $T_mM$. Given a subgroup $K$ of $G$, there is a (left) $K$-action on the product $G\times N$ given by 
\begin{equation}
	k\cdot (g,\nu)=(gk^{-1},k\cdot \nu).
\end{equation}
This action is free and proper by freeness and properness of the action on the first factor. The orbit space $G\times_KN$ is thus a smooth manifold whose points are equivalence classes of the form $[(g,\nu)]$, and the orbit map $\rho:G\times N\to G\times_KN$ is a smooth surjective submersion. Moreover the group $G$ acts smoothly and properly on $G\times_KN$, by left multiplication on the first factor. The theorem below states that, given a $G_m$-invariant neighbourhood of zero $N_0\subset N$, the associated bundle $G\times_KN_0$ defines a \defn{local model} for some $G$-invariant open subset $U\subset M$.

\begin{theorem}[Tube Theorem (cf. \cite{MR2021152} Theorem $2.3.28$)]\label{Tube theorem}
	Let $(M,G)$ be a proper $G$-manifold and set $K=G_m$ for some $m\in M$.  Let $N_0\subset N$ be an open $K$-invariant neighbourhood of $0$. Then, there exists a $G$-invariant neighbourhood $U\subset M$ of $m$ and a $G$-equivariant diffeomorphism
	\begin{equation}\label{tau tube}
		\varphi:G\times_K N_0\to U
	\end{equation}
sending $[(e,0)]$ on $m$.
\end{theorem}

\subsection*{G-categorical open subsets}
Although the $LS$-category is defined for topological spaces, our interest is mainly about proper $G$-manifolds $(M,G)$. The terminologies are thus defined in this setting. Given $(M,G)$, a homotopy $H:M\times [0,1]\to M$ which satisfies $H(g\cdot m,t)=g\cdot H(m,t)$ for every $g\in G$, $m\in M$ and $t\in [0,1]$ is called a \defn{$G$-homotopy}. We write $H_t(m)=H(m,t)$. Let $B\subset A$ be two $G$-invariant subsets of $M$. A \defn{$G$-deformation retract of $A$ onto $B$} is a $G$-homotopy $H:A\times [0,1]\to A$ such that $H_0(a)=a$ and $H_1(a)\in B$ for every $a\in A$, and $H(b,1)=b$ for every $b\in B$.

\begin{definition}\label{def: G-categorical open subset}
\normalfont A $G$-invariant subset $U\subset M$ is called \defn{$G$-categorical} if there exists a $G$-deformation retract of $U$ onto the orbit $G\cdot x$ of some $x\in U$.
\end{definition}

\begin{definition}
	\normalfont The \defn{equivariant LS-category of $M$}, denoted $\cat_G(M)$, is the least number of $G$-categorical open subsets $U\subset M$ that are required to cover $M$. We set $\cat_G(M)=\infty$ if such a cover does not exist. The non-equivariant LS-category $\cat(M)$ is obtained by setting $G=\mathbbm{1}$.
\end{definition}

\begin{example}\label{ex:TETRAHEDRON}
\normalfont The equivariant version of the LS-category is in general different from its non-equivariant analogue, as shown in the examples below.
	\begin{enumerate}[label=(\roman*)]
		\item Let $M=S^1\times \R$ with cylindrical coordinates $(\theta,z)$. Define an $S^1$-action on $M$ by $\phi\cdot (\theta,z)=(\theta+\phi,z)$. The cylinder itself is an $S^1$-categorical open subset with $S^1$-deformation retract $H:M\times [0,1]\to M$ given by $H((\theta,z),t)=(\theta,(1-t)z)$. Therefore, $\cat_{S^1}(M)=1$. However we require two contractible open subsets to cover $M$, which yields $$1=\cat_{S^1}(M)<\cat(M)=2.$$

                                                                                                                                                             		\item Consider the complex projective space $M=\C P^2$ endowed with the $S^1$-action $\theta\cdot [z_0:z_1:z_2]=[e^{i\theta}z_0:z_1:z_2].$ For $i=0,1,2$, the open subsets $U_i=\lbrace [z_0:z_1:z_2]\mid z_i\neq 0\rbrace$ are $S^1$-invariant. On $U_0$, an $S^1$-deformation retract onto an orbit is given by
                                                                                                                                                  $$H([z_0:z_1:z_2],t)=[z_0:(1-t)z_1:(1-t)z_2].$$ The image $H_1(U_0)$ is the single point $[1:0:0]$ which is a fixed point of the action, hence an $S^1$-orbit. Similar homotopies can be found on $U_1$ and $U_2$, respectively. Therefore $\cat_{S^1}(M)$ is at most three. The fact that we have an equality follows from Proposition \ref{atmost} below. We conclude that $$\cat_{S^1}(M)=\cat(M)=3.$$
\item 
The rotations of a tetrahedron form a group $\mathbb{T}$ of order $12$, which is a zero-dimensional Lie subgroup of $SO(3)$. This group acts on $M=S^2$. We construct a cover of $M$ by three $\mathbb{T}$-categorical open subsets as follows:

Pick a point $x_1\in M$ and its opposite point $y_1\in M$. The $\mathbb{T}$-orbit of $x_1$ forms a spherical tetrahedron with vertices $x_1,x_2,x_3,x_4$. Similarly the $\mathbb{T}$-orbit of $y_1$ forms another spherical tetrahedron with vertices $y_1,y_2,y_3,y_4$. For each $i<j$ denote by $p_{ij}$ the middle point of the geodesic arc joining $x_i$ and $x_j$ (cf. Figure \ref{spherical tetrahedrons}).

\begin{figure}[!ht]
	\centering	
		\begin{minipage}[t]{9cm}
		\centering
			\includegraphics[height=5cm]{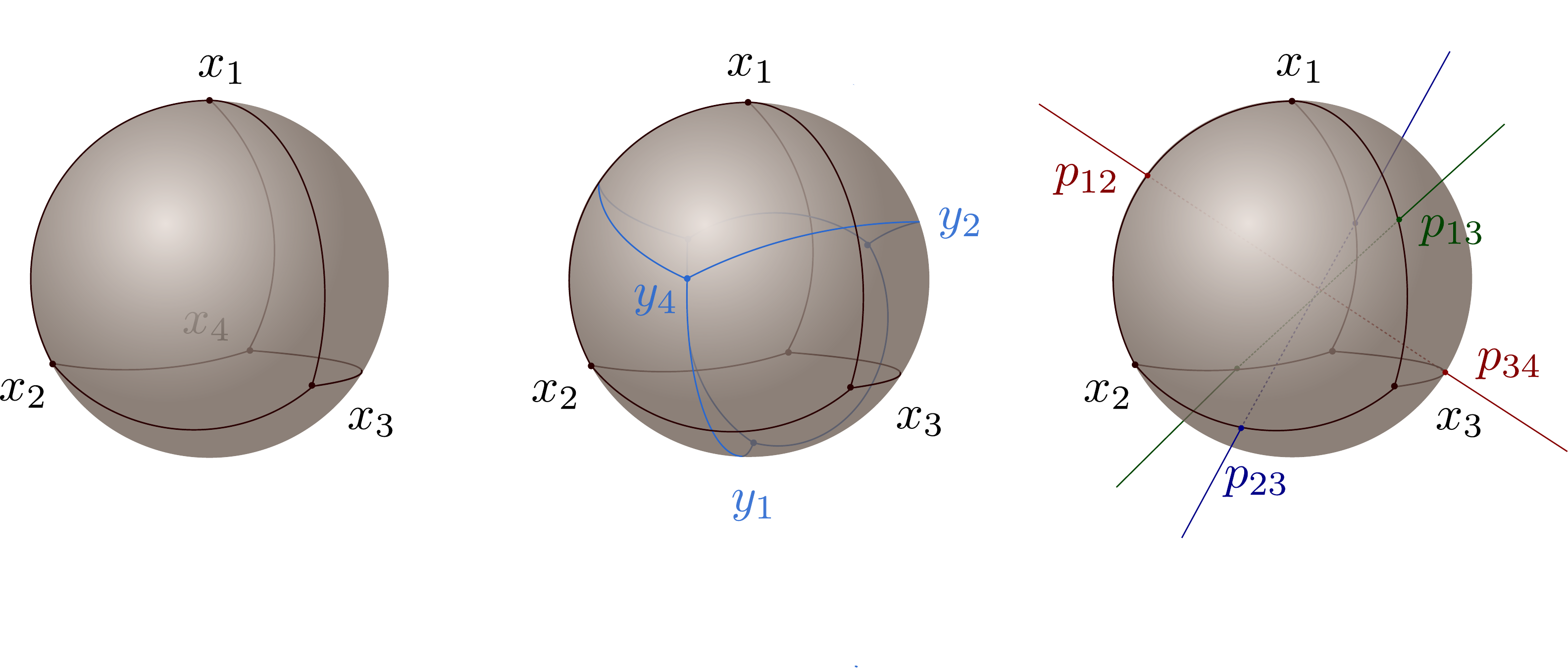}
	\caption{\small Spherical tetrahedrons on the sphere.}
	\label{spherical tetrahedrons}
\end{minipage}
\end{figure}

For each $i$, let $D_i\subset M$ be an open disk centered at $x_i$ such that
 $$\overline{D_i}\cap\overline{D_j}=\lbrace p_{ij}\rbrace\quad \forall i<j. $$
 In the same way, let for each $i$, an open disk $E_i\subset M$ centered at $y_i$ with the property
 $$\forall i<j\quad \overline{E_i}\cap\overline{E_j}=\lbrace p_{kl}\rbrace\quad \mbox{where} \quad k,l\notin \lbrace i,j\rbrace\quad k<l $$ as shown in Figure \ref{spherical}.

\begin{figure}[!ht]
	\centering	
		\begin{minipage}[t]{9cm}
		\centering
		\includegraphics[height=4cm]{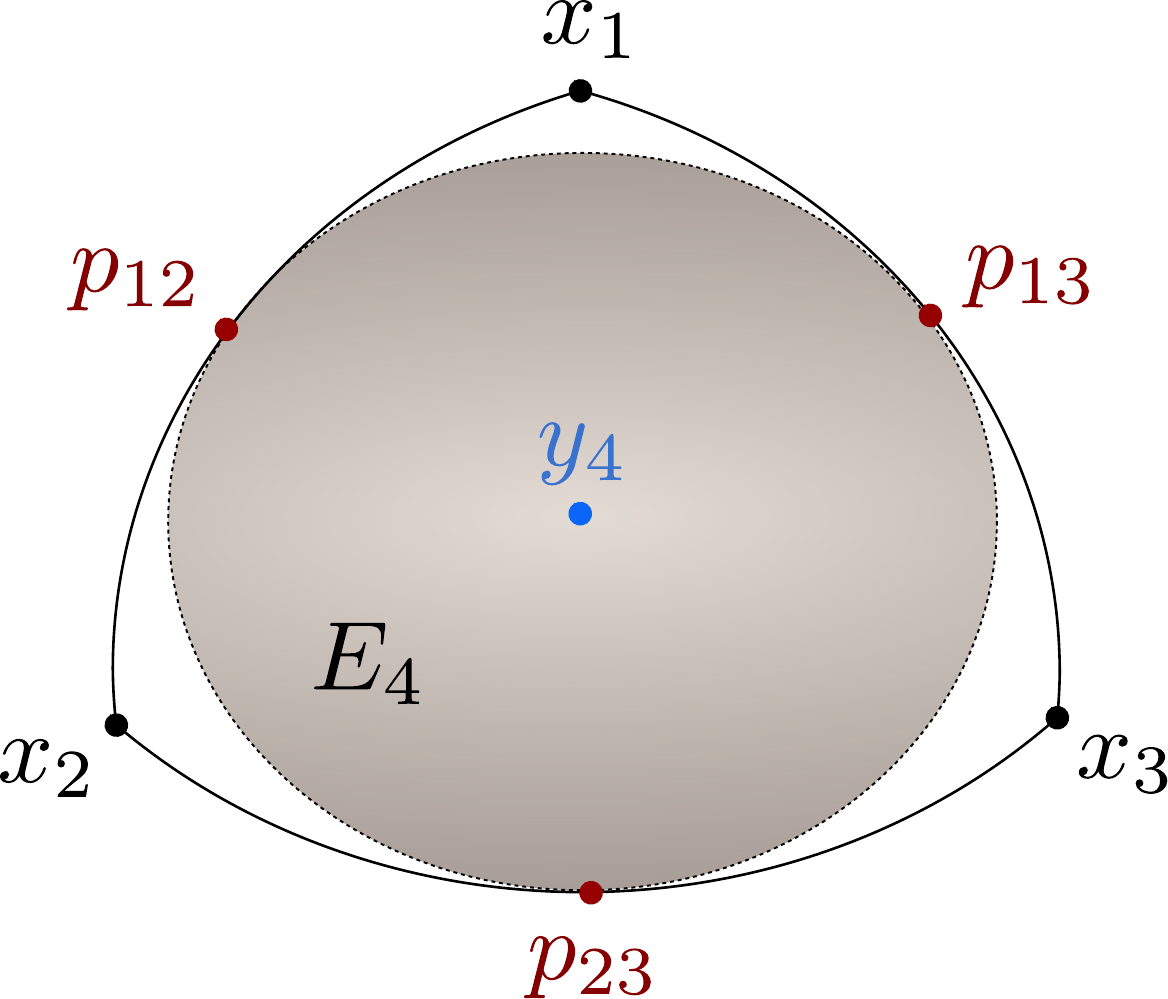}	
		\caption{\small Disk $E_4$ centered at $y_4$.}
		\label{spherical} 
\end{minipage}	
\end{figure}

 Finally we define for each $i<j$, an open subset $B_{ij}\subset M$ containing $p_{ij}$ such that
 
 \begin{eqnarray*}
 x_k\in\overline{B_{ij}}\setminus B_{ij} && \forall k=i,j\\
 y_k\in\overline{B_{ij}}\setminus B_{ij} && \forall k\neq i,j
 \end{eqnarray*}

We obtain the following  $\mathbb{T}$-categorical open subsets:
$$D=\bigcup\limits_{i=1}^4 D_i$$ which retracts in a $\mathbb{T}$-equivariant way onto the orbit $\mathbb{T}\cdot x_1,$
$$E=\bigcup\limits_{i=1}^4 E_i$$ which retracts in a $\mathbb{T}$-equivariant way onto the orbit $\mathbb{T}\cdot y_1,$ and
$$B=\bigcup\limits_{i<j} B_{ij}$$ which retracts in a $\mathbb{T}$-equivariant way onto the orbit $\mathbb{T}\cdot p_{12}.$ Those three subsets form a cover of $M$. This cover is in fact the smallest that we can take, by Proposition \ref{atmost} below. Hence $$3=\cat_{\mathbb{T}}(M)>\cat(M)=2.$$ \label{tetrahedron}
\end{enumerate}
\end{example}

\section{Stratifications and orbit-type strata} \label{isotropytype}
A \defn{partition} of a topological space $M$ is a cover of $M$ by pairwise disjoint subsets. Clearly every topological space admits a partition into its connected components. If our topological space is endowed with a group action, we can choose a partition which also encodes the information about the group action. For example, a proper $G$-manifold $(M,G)$ can be partitioned into locally closed (connected) submanifolds called the orbit-type strata, each of them being a union of group orbits with the same orbit-type.

\subsection*{Stratifications}\label{s: stratification}
There are several ways to define stratifications. The one we present here is the definition used by\cit{Kirwan} in her thesis \cite{Kirwan}. It is more flexible than the standard definition of\cit{Duistermaat and Kolk} \cite{kolk} (Definition $2.7.3$), especially for applications to algebraic geometry. A \defn{strict partial order} $\prec$ on $\mathcal{B}$ is a binary relation which is irreflexive and transitive. Note that in this case if $\alpha,\beta\in\mathcal{B}$ are such that $\alpha\prec\beta$, then $\beta\not\prec\alpha$.
For example the set of conjugacy classes of subgroups of $G$ admits the strict partial order $\prec_{conj}$, where we say that $(K)\prec_{conj} (H)$ if and only if $H$ is conjugate to a proper subgroup of $K$.

\begin{example}
\normalfont The group $\mathbb{T}$ has four conjugacy classes, namely $(\mathbb{T}),(\Z_3),(\Z_2)$ and $(\mathbbm{1})$. There are partially ordered with respect to $\prec_{conj}$ as shown in Figure \ref{T}.

\vspace{0.5cm}
\centering
\begin{figure}[h!]
\centering
\begin{minipage}[t]{9cm}
	\[
			\begin{xy}
					(0,0)*+{\mathbb{T}}="T";%
			 		(-10,-10)*+{\Z_3}="z";(10,-10)*+{\Z_2}="zz";%
			 		(0,-20)*+{\mathbbm{1}}="i";%
					{\ar@{-} "T";"zz"};{\ar@{-} "T";"z"};%
					{\ar@{-} "zz";"i"};{\ar@{-} "z";"i"};%
			\end{xy}
	\]\caption{\small Conjugacy classes of subgroups of $\mathbb{T}$ where the order goes up to down i.e. $(\mathbb{T})$ is minimal with respect to $\prec_{conj}$.}
	\label{T}
	\end{minipage}
	\end{figure}
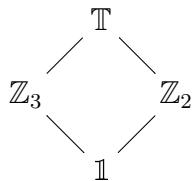
\end{example}

\begin{definition}
			\normalfont A collection $\left\lbrace M_{\beta}\mid \beta\in \mathcal{B}\right\rbrace$ of subsets of a topological space $M$ is \defn{locally finite} if each compact subset of $M$ meets only finitely many $M_{\beta}$. A locally finite collection $\left\lbrace M_{\beta}\mid \beta\in \mathcal{B}\right\rbrace$ of locally closed (non-empty) topological subspaces of $M$ form a \defn{$\mathcal{B}$-decomposition} or \defn{stratification} of $M$ if $M$ is the disjoint union of the \defn{strata} $M_{\beta}$, and there is a strict partial order $\prec$ on the indexing set $\mathcal{B}$ such that
\begin{equation}\label{closure condition}
	\overline{M}_{\beta}\subset\bigcup_{\alpha\preceq\beta}M_{\alpha}
\end{equation}
for every $\beta\in\mathcal{B}$. We say that the stratification is \defn{smooth} if $M$ is a smooth manifold and every $M_{\beta}$ is a locally closed submanifold. 
\end{definition}
Given a stratification of $M$, a strict partial order can be defined on the strata as follows
\begin{equation}\label{partial order strata}
M_{\alpha}< M_{\beta} \quad \Longleftrightarrow\quad \alpha\prec\beta.
\end{equation}
We say that a stratum $M_{\beta}$ is \defn{minimal} with respect to \eqref{partial order strata} if there is no $\alpha\in\mathcal{B}$ such that $M_{\alpha}< M_{\beta}$. Of course minimal strata are not unique because we just have a partial ordering. If $M_{\beta}$ is a minimal stratum, \eqref{closure condition} implies that $\overline{M}_{\beta}\subset M_{\beta}$. In particular $\overline{M}_{\beta}=M_{\beta}$ i.e. $M_{\beta}$ is closed in $M$.
\subsection*{Orbit-type strata}
Let $G$ be a Lie group and $H\subset G$ be a closed subgroup. The \defn{conjugacy class} of $H$ is the set $(H)=\lbrace L\subset G\mid L=gHg^{-1}\mbox{ for some }g\in G\rbrace.$ 
Given a $G$-manifold $(M,G)$, we define the set $$M_{(H)}:=\left\lbrace m\in M\mid G_m\in (H)\right\rbrace$$ which is the union of all the $G$-orbits in $M$ with orbit-type $(H)$. Using the definitions and the $G$-invariance of $M_{(H)}$, it is shown in \cite{MR2021152} (Proposition $2.4.4$) that $M_{(H)}=G\cdot M_H$ where $$M_H=\lbrace m\in M\mid G_m=H\rbrace.$$ Note that the biggest subgroup of $G$ which leaves $M_H$ invariant is the normalizer $N_G(H)=\left\lbrace g\in G\mid gHg^{-1}=H\right\rbrace.$ Furthermore this action induces a well-defined free action of the quotient group $N_G(H)/H$ on $M_H$. Write
$$M_H=\coprod_{b\in\mathcal{B}_H} M_{H,b}$$ as the disjoint union of its connected components, indexed on some set $\mathcal{B}_H$. Given $b\in \mathcal{B}_H$, we define the equivalence class $(b)$ to be the set of indices $a\in \mathcal{B}_H$ such that $G\cdot M_{H,a}=G\cdot M_{H,b}$. Let $\mathcal{B}$ be the set of pairs $\beta=((H),(b))$ where $(H)$ is the conjugacy class of some closed subgroup of $G$ and $b\in \mathcal{B}_H$.

\begin{definition}\label{def: orbit type stratum}
\normalfont For $\beta=((H),(b))\in \mathcal{B}$, we define an \defn{orbit-type stratum} $M_{\beta}$ to be the $G$-orbit of the connected component $M_{H,b}$ of $M_H$. 
\end{definition}
We use here a modified definition of the standard definition which states that an orbit-type stratum is a connected component of $M_{(H)}$. If the $G$-action on $M$ is proper, the connected components $M_{H,b}$ are locally closed embedded submanifolds of $M$ and so are their $G$-orbits (cf. \cite{MR2021152} Proposition $2.4.7$).

The example below illustrates the difference between the standard definition of orbit-type strata and ours. With our definition, the orbit-type strata might not be connected.
\begin{example}
	\normalfont Think of $\R^*=\R\setminus\lbrace 0\rbrace$ as a multiplicative group and let it act on $M=\R^2$ by $t\cdot (x,y)=(x,ty)$. The stabilizers of points of $M$ are either equal to the trivial group $\mathbbm{1}$, or equal to $\R^*$. Then $M_{(\R^*)}=M_{\R^*}$ is the $x$-axis, and $M_{(\mathbbm{1})}=M_{\mathbbm{1}}=H_+\cup H_-$ where $H_{\pm}=\lbrace (x,y)\in M\mid \pm y > 0\rbrace$. According to the standard definition of orbit-type strata, there are two strata with orbit-type $(\mathbbm{1})$, namely the connected components $H_+$ and $H_-$; and one stratum with orbit-type $(\R^*)$, the $x$-axis.

With our definition, there is one stratum with orbit-type $(\R^*)$ which is the $x$-axis; but there is only one stratum with orbit-type $(\mathbbm{1})$ which is $H_+\cup H_-$. Indeed, $M_{\mathbbm{1}}$ has two connected components, $H_+$ and $H_-$. The $\R^*$-orbits of each of them coincide. There is thus only one stratum with orbit-type $(\mathbbm{1})$ and it is not connected.
\end{example}
 We define a strict partial order $\prec$ on $\mathcal{B}$ as follows: for $\alpha=((K),(a))$ and $\beta=((H),(b))$,
\begin{equation}\label{partial order orbit}
	\alpha\prec \beta\quad\Longleftrightarrow\quad \alpha\neq\beta\quad\mbox{and}\quad M_{\alpha}\cap\overline{M}_{\beta}\neq \varnothing.
\end{equation}
 By $\alpha\neq\beta$ we mean that the associated orbit-type strata $M_{\alpha}$ and $M_{\beta}$ are distinct. 
 
\begin{proposition}[\cit{Sjamaar and Lerman} \cite{MR1127479}]\label{stratification by orbit type}
		Let $(M,G)$ be a proper $G$-manifold and let $(\mathcal{B},\prec)$ as above with partial order \eqref{partial order orbit}. Then the orbit-type strata $\lbrace M_{\beta}\mid \beta\in\mathcal{B}\rbrace$ form a smooth stratification of $M$.
\end{proposition}
\begin{remark}	In \cite{MR1127479}, \cit{Sjamaar and Lerman} use the standard definition of orbit-type strata. In particular those are connected submanifolds of $M$. In \cite{moi} we give a proof using the modified definition of orbit-type strata (cf. \cite{moi} Proposition $2.2.9$).
\end{remark}

\begin{example}\label{ex: stratum}\normalfont
	Given an equivalence class $(H)$, the corresponding orbit-type strata might not all have the same dimension, as shown in the following example, appearing in\cit{Delzant} \cite{delzant} and\cit{Sjamaar and Lerman}~\cite{MR1127479}.
	\begin{enumerate}[label=(\roman*)]
	\item  Let $M=\C P^2$ endowed with the $S^1$-action $$\theta\cdot [z_0:z_1:z_2]=[e^{i\theta}z_0:z_1:z_2].$$ The set $M_{S^1}$ has two connected components namely, the point $[1:0:0]$ and a copy of $\C P^1$, which consists of the points of the form $[0:z_1:z_2]$. Since $S^1$ acts trivially on each of these components, they form themselves two orbit-type strata, which are closed submanifolds of $M$. Since the action is free anywhere else, the last orbit-type stratum is $M\setminus \left(\lbrace [1:0:0]\rbrace \cup \C P^1\right)$. It has orbit-type $(\mathbbm{1})$ and is an open dense submanifold of $M$. 
	
	\item Let $M=S^2$ be the $2$-sphere embedded in $\R^3$, equipped with the $S^1$-action which rotates the sphere about the $z$-axis. There are three orbit-type strata namely, $M_{(\mathbbm{1})}$ which is diffeomorphic to $S^1\times (-1,1)$ and the two closed connected components of $M_{S^1}$ that are the North and South pole. \label{ex: stratum sphere}
	\item \label{ex: stratum T}  The group $\mathbb{T}$ acts on $M=S^2$. This group contains a copy of the cyclic group of order three $C_3\simeq \Z_3$ for each vertex, one copy of $\Z_2$ for each axis joining the middle point of an edge and the middle point of the opposite edge, and the identity (cf. Figure ~\ref{fig:TT}).
	
	\vspace{0.3cm}
	\begin{figure}[!ht]
		\centering	
		\begin{minipage}[t]{9.5cm}
			\centering
			\includegraphics[height=4cm]{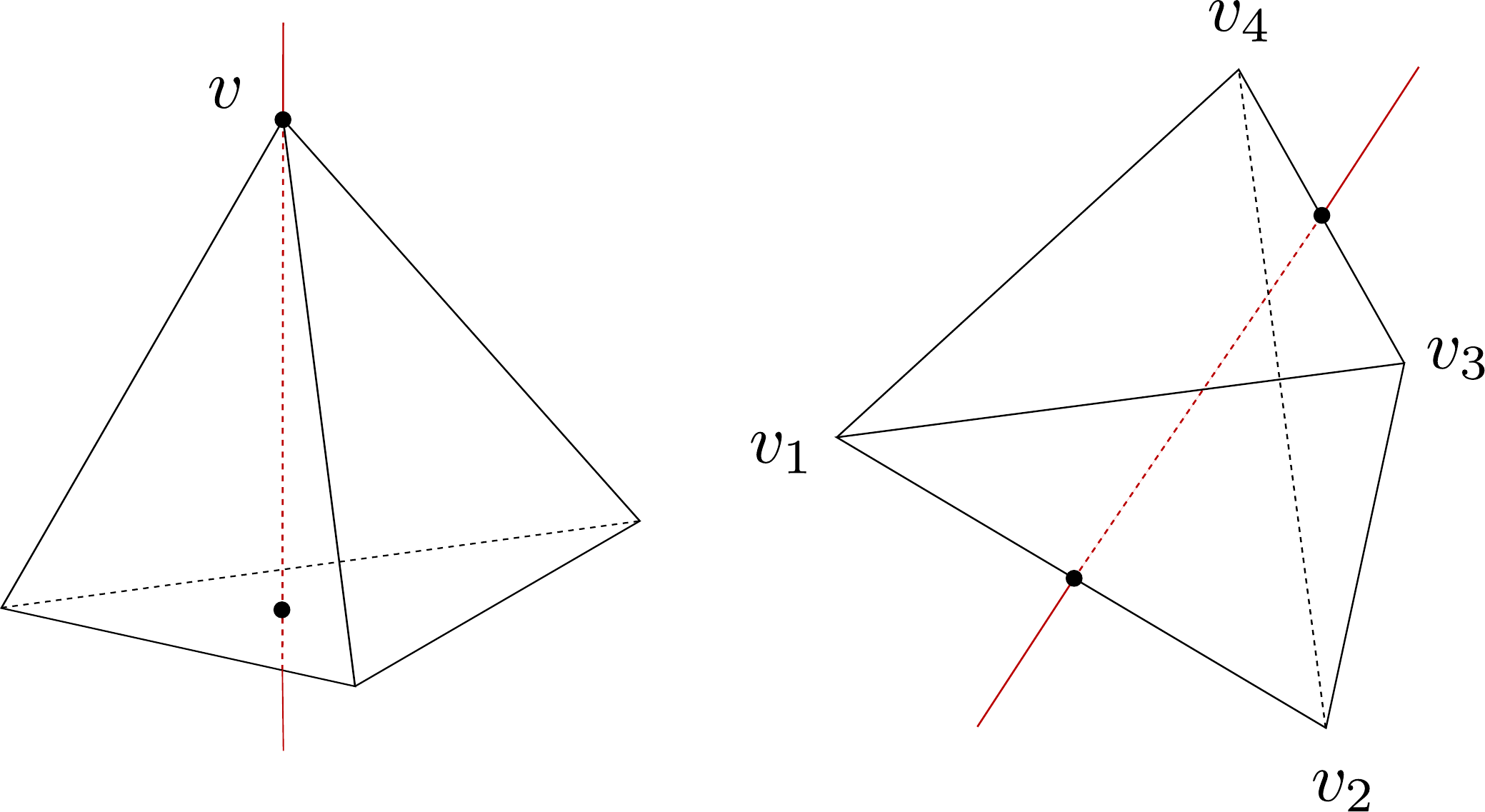}
			\caption{\small On the left hand side we fix a vertex $v$ and permute the three other vertices. As a subgroup it is isomorphic to $C_3$. On the right hand side we permute $v_1,v_2$ and $v_3,v_4$. This subgroup is isomorphic to $\Z_2$.}
			\label{fig:TT}
		\end{minipage}
	\end{figure}

There are two minimal strata with orbit-type $(\Z_3)$, one minimal stratum with orbit-type $(\Z_2)$, and one open dense stratum with orbit type $(\mathbbm{1})$ (cf. Figure \ref{ex: stratum tetrahedron}). Indeed, when $H=\Z_3$, the eight points forming $M_{(H)}$ are a union of two $\mathbb{T}$-orbits. There are thus two strata with orbit-type $(\Z_3)$. For $H=\Z_2$, the six points forming $M_{(H)}$ are a single $\mathbb{T}$-orbit and form a single stratum.

\vspace{0.3cm}
\begin{figure}[h!]
\begin{center}	
	\begin{tabular}{|c|c|c|}
  			\hline
  			$H$ & $M_H$  & $M_{(H)}$  \\
  			\hline
  			$\mathbb{T}$ & $\varnothing$ & $\varnothing$ \\
  			\hline
  			&&\\
  			\raisebox{25pt}{$\Z_3$} & \includegraphics[height=2cm]{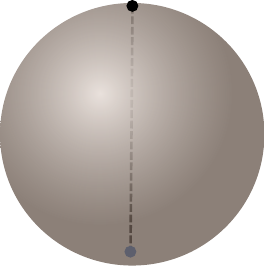} & \includegraphics[height=2cm]{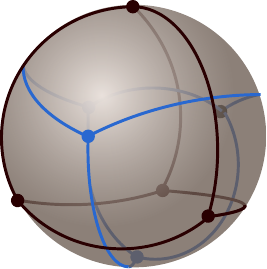}\\
  			&&\\
  			\hline
  			&&\\
  			\raisebox{30pt}{$\Z_2$}& \includegraphics[height=2cm]{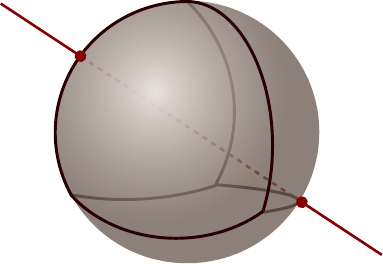} & \includegraphics[height=2.5cm]{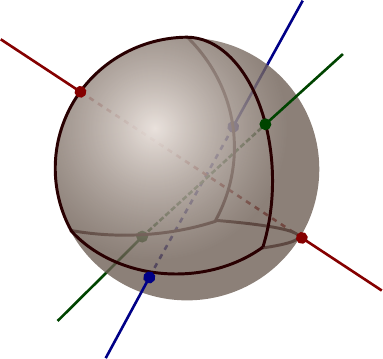}\\
  			&&\\
  			
  			\hline
  			$\mathbbm{1}$& $M\setminus\lbrace 14\mbox{ points }\rbrace$& $M\setminus\lbrace 14\mbox{ points }\rbrace$\\
  			\hline
	\end{tabular}
	\caption{\small Orbit-type strata of $(M,\mathbbm{T})$ where $M=S^2$.}\label{ex: stratum tetrahedron}
\end{center}
\end{figure}

	\end{enumerate}
\end{example}

In general, using the strict partial order $\prec_{conj}$ on the conjugacy classes of subgroups of $G$ is not enough to guarantee that we have a good stratification. For instance, in Example \ref{ex: stratum} \ref{ex: stratum T}, we have $(\mathbb{T})\prec_{conj}(\Z_2)$ but there are no strata with orbit-type $(\mathbb{T})$. However we have the following lemma:

\begin{lemma}\label{lemma conj}
	If $\alpha=((K),(a))$ and $\beta=((H),(b))$ then $$\alpha\prec\beta\quad\Longrightarrow\quad (K)\prec_{conj} (H).$$
\end{lemma}
\begin{Proof}
	By definition $\alpha\prec\beta$ implies that there exists some $x\in M_{\alpha}\cap\overline{M}_{\beta}$. In particular $x\in M_{\alpha}$ and then $G_x\in (K)$. By the Tube Theorem \ref{Tube theorem}, there is a $G$-invariant open neighbourhood $U\subset M$ of $x$, locally modelled by a fixed associated bundle $G\times_{G_x}N_0$, in which $x$ reads $[(e,0)]$.
	
	By definition of the adherence, there is a sequence $(x_n)_{n\in\N}\subset M_{\beta}$ converging to $x$ in $M$, with stabilizers $G_{x_n}\in (H)$. For $n$ big enough, $x_n\in U$ and it can thus be identified with some point $[(g_n,\nu_n)]\in G\times_{G_x}N_0$. The stabilizer of $[(g_n,\nu_n)]$ is $$G_{[(g_n,\nu_n)]}=g_n(G_x)_{{\nu}_n}g_n^{-1}$$ and is thus conjugate to a proper subgroup of $G_x$, because by assumption $M_{\alpha}$ and $M_{\beta}$ are disjoint. Since $G_x\in (K)$ and $G_{x_n}\in (H)$, it follows that $(K)\prec_{conj}(H)$.
\end{Proof}

\subsection*{Stratification of a convex polytope by open faces}\label{stratification polytope}
There is a natural stratification of a convex polytope into vertices, edges and higher dimensional faces. Let $\Delta\subset (\R^n)^*$ be a $n$-dimensional convex polytope. Let $X_1,\dots,X_d$ in $\R^n$ be the outward-pointing normal vectors to the facets. Then there exists real numbers $\lambda_1,\dots,\lambda_d$ such that $\Delta$ reads
\begin{equation*}
	\Delta=\bigcap_{i=1}^d\left\lbrace \mu\in (\R^n)^*\mid \langle \mu,X_i\rangle\leq\lambda_i\right\rbrace.
\end{equation*}

Let $\mathcal{B}$ be the set of subsets (possibly empty) $\beta\subset\lbrace 1,\dots,d\rbrace.$ For each $\beta\in\mathcal{B}$ we consider the intersection
\begin{equation*}
	F_{\beta}=\bigcap_{i\in \beta}\left\lbrace \mu\in\Delta\mid \langle \mu,X_i\rangle=\lambda_i\right\rbrace.
\end{equation*}
If $F_{\beta}\neq\varnothing$, its relative interior $\overset{\circ}{F}_{\beta}$ is called a $l$-dimensional open face of $\Delta$ where $l$ is equal to $n$ minus the cardinality of $\beta$. We equip $\mathcal{B}$ with the strict partial order
\begin{equation}\label{partial order polytope}
	\alpha\prec\beta\quad\Longleftrightarrow\quad \alpha\neq\beta\quad\mbox{and}\quad \overset{\circ}{F}_{\alpha}\cap \overline{\overset{\circ}{F}_{\beta}}\neq\varnothing.
\end{equation}
With this strict partial order, the collection $\lbrace \overset{\circ}{F}_{\beta}\mid \beta\in\mathcal{B}\rbrace$ forms a $\mathcal{B}$-stratification of $\Delta$. A strict partial order is defined on the set of faces by
\begin{equation*}
	\overset{\circ}{F}_{\alpha}< \overset{\circ}{F}_{\beta}\quad \Longleftrightarrow\quad\alpha\prec\beta.
\end{equation*}
Finally note that, if $\alpha\prec \beta$ then $\beta\subset \alpha$.

\section{G-tubular covers}\label{s:G-categorical tubular covers}
If $(M,G)$ is a proper $G$-manifold, the Tube Theorem \ref{Tube theorem} allows us to produce $G$-categorical open subsets in the following way: any $m\in M$ admits a $G$-invariant neighbourhood $U\subset M$ such that the map $\varphi:Y_0\to U$ defined in \eqref{tau tube} is a $G$-equivariant diffeomorphism. Here $$Y_0=G\times_{G_m}N_0$$ where $N_0$ is a fixed neighbourhood of zero in some subspace $N\subset T_mM$, complementary to $\mathfrak{g}\cdot m$ in $T_mM$, on which $G_m$ acts linearly. The proper $G$-manifold $Y_0$ is a local model for $U$, in which $m$ reads $\varphi^{-1}(m)=[e,0]$. The $G$-homotopy $F:Y_0\times [0,1]\to Y_0$ defined by
\begin{equation*}
	F([(g,\nu)],t)=[(g,(1-t)\nu)].
\end{equation*}
is a $G$-deformation retract of $Y_0$ onto the orbit $G\cdot [e,0]$. By using the fact that $\varphi$ is a $G$-equivariant diffeomorphism, the open subset $U=\varphi(Y_0)$ is $G$-categorical since the $G$-homotopy $H:U\times [0,1]\to U$ given by 
\begin{equation}\label{Gtubular homotopy}
	H(p,t)=\varphi\left(F(\varphi^{-1}(p),t)\right)
\end{equation}
is a $G$-deformation retract of $U$ onto $G\cdot m$.
\begin{definition}\label{def:G tubular open subset}
	\normalfont A $G$-categorical open subset $U\subset M$ as above, with associated $G$-deformation retract as in \eqref{Gtubular homotopy}, is called a \defn{$G$-tubular open subset} of $M$. A cover of $M$ made of $G$-tubular open subsets is called a \defn{$G$-tubular cover} of $M$.
\end{definition}

\subsection*{Minimal G-tubular covers}
Clearly, every $m\in M$ admits a neighbourhood which is a $G$-tubular open subset of $M$. Consequently, $G$-tubular covers of $M$ always exist. The question is whether they can be refined. Let $\mathcal{U}$ be any $G$-tubular cover of $M$. We know that $M$ can be decomposed into the disjoint union of its orbit-type strata $\lbrace M_{\beta}\mid \beta\in\mathcal{B}\rbrace$, which form themselves a smooth stratification of $M$. Let $\mathcal{B}'\subset\mathcal{B}$ be the biggest subset of indices $\beta\in\mathcal{B}$ such that $M_{\beta}$ is minimal with respect to \eqref{partial order strata}. Consider the disjoint union $M'$ of all the strata $M_{\beta}$ with $\beta\in\mathcal{B}'$. From $\mathcal{U}$ we extract a subcover $\mathcal{U}'$, chosen as small as possible such that $\mathcal{U}'$ covers $M'$. In particular $\mathcal{U}'$ is a refinement of $\mathcal{U}$. We ask the following:

\vspace{0.5cm}
\hypertarget{label}{(Q)} \textit{Does it exist a subset $\mathcal{U}'\subset \mathcal{U}$, obtained as above, which covers $M$?} 
\vspace{0.5cm}

The answer is in general negative (cf. Example \ref{non-example}). However it is positive for all the proper $G$-manifolds listed in Example \ref{ex:TETRAHEDRON}, where $\mathcal{U}'$ is constructed explicitly.
\begin{definition}\label{minimal tubular cover}
\normalfont Let $(M,G)$ be a proper $G$-manifold. $\mathcal{U}'$ obtained as above is called a \defn{minimal $G$-tubular cover} if the following holds:
\begin{enumerate}[label=(\roman*)]
	\item $\mathcal{U}'$ is a cover of $M$.\label{minimal i}
	\item For each minimal orbit-type stratum $M_{\beta}$, the set $$\mathcal{V}'_{\beta}=\lbrace V_{\beta}=U\cap M_{\beta}\mid U\in\mathcal{U}'\rbrace$$ is the smallest cover by $G$-categorical open subsets of $M_{\beta}$, where the topology of $M_{\beta}$ is the subset topology.\label{minimal ii}
\end{enumerate}	
\end{definition}
We discuss the simplest example where such a cover exists. Let $S^2\subset\R^3$, on which $S^1$ acts by rotations about the $z$-axis. This action has two minimal orbit-type strata, namely the North and South pole. Two small disks centered at those points are $S^1$-tubular open subsets and can be taken sufficiently big so that they form a minimal $S^1$-tubular cover of $S^2$. In this example, a disk centered at the North pole can be extended until its closure meets the South pole. The impossibility to extend it further relies on the fact that such neighbourhoods are constructed by mean of the Riemannian exponential map. This map is no longer injective if the disk contains two opposite points on the sphere. The next proposition gives another answer to this fact by using the properties of $G$-tubular open subsets.

\begin{proposition}\label{atmost}
		Let $(M,G)$ be a proper $G$-manifold. If $U\subset M$ is a $G$-tubular open subset which intersects a minimal orbit-type stratum $M_{\beta}$, then $U$ retracts onto the orbit $G\cdot x$ of some $x\in M_{\beta}$. In particular $G$-tubular open subsets intersect at most one minimal orbit-type stratum.
\end{proposition}
\begin{Proof}
	Let $\beta=((H),(b))\in\mathcal{B}$ such that $M_{\beta}$ is a minimal orbit-type stratum. Let $U\subset M$ be a $G$-tubular open subset of $M$ such that $U\cap M_{\beta}\neq\varnothing$, and let $H:U\times [0,1]\to U$ be a $G$-deformation retract of $U$ onto $G\cdot x$ for some $x\in M$. By contradiction, assume that $x\in M_{\alpha}$ for some $\alpha=((G_x),(a))\neq\beta$. 
	
 Each point $y\in U\cap M_{\beta}$ has stabilizer $G_y\in (H)$. By $G$-equivariance of the homotopy, $G_y$ is a subgroup of $G_{H_1(y)}$ which is itself conjugate to $G_x$, as $H_1(y)$ and $x$ lie on the same orbit. In particular $(G_x)\prec_{conj}(H)$. Two cases occur:
 \begin{enumerate}[label=(\roman*)]
 	\item If $M_{\beta}\cap\overline{M}_{\alpha}\neq\varnothing$, then $\beta\prec\alpha$ since $\beta\neq\alpha$. By Lemma \ref{lemma conj} we get $(H)\prec_{conj}(G_x)$ which is a contradiction.
 	\item If $M_{\beta}\cap\overline{M}_{\alpha}=\varnothing$ we must use the assumption that $U$ is $G$-tubular. Let $G\times_{G_x}N_0$ be the local model for $U$. Given $y\in M_{\beta}$ we define the $G$-equivariant path $y(t)=H_t(y)$, where $t\in [0,1]$. In the local model, $y$ reads $[g,\nu]$ and $y(t)$ reads $[g,\nu_t]$ where $\nu_t=(1-t)\nu$. We can assume without lost of generality that $(G_x)_{\nu}=H$. Observe that, by linearity of the $G_x$-action on $N_0$, we have $(G_x)_{\nu_t}=(G_x)_{\nu}=H$ for all $t\neq 1$. Hence $G_{[g,\nu_t]}=g(G_x)_{\nu_t}g^{-1}=gHg^{-1}$ for every $t\neq 1.$ In particular, $y(t)\in M_{\beta}$ for all $t\neq 1$.
 	Since the path $y(t)$ starts at $y\in M_{\beta}$ and ends on $G\cdot x\subset M_{\alpha}$, there is some $t_0\in [0,1]$ such that $y(t_0)\in \overline{M_{\alpha}}$. The parameter $t_0$ is chosen the smallest such that this occurs. If $t_0\neq 1$, the previous argument shows that $y(t_0)\in M_{\beta}\cap\overline{M_{\alpha}}$, which is a contradiction. Otherwise, since $y(t)\in M_{\beta}$ for all $t<t_0$, there is a sequence $(y_n)_{n\in\N}\subset M_{\beta}$ which converges to $y(t_0)$. By closedness of $M_{\beta}$, this yields $y(t_0)\in M_{\beta}\cap \overline{M_{\alpha}}$, which is again a contradiction.  We conclude that $x\in M_{\beta}$. 
 \end{enumerate}
\end{Proof}

The answer to question \hyperlink{label}{(Q)} is in general negative. In the example below, $(M,G)$ is a proper $G$-manifold with $M$ compact, and the action admits only one minimal orbit-type stratum $M_{\beta}$ which is a single $G$-orbit.

\begin{example}\label{non-example}
	\normalfont 
	
Think of $M=S^3$ as the set of unit vectors $(z_1,z_2)\in \C^2$ equipped with the $S^1$-action
\begin{equation*}
		\theta\cdot (z_1,z_2)=(e^{i\theta}z_1,e^{2i\theta}z_2).
\end{equation*}
	This action has only one minimal orbit-type stratum $M_{\beta}$ with $\beta=((\Z_2),(b))$ for some index $b\in \mathcal{B}_{\Z_2}$. Explicitly
	\begin{equation*}
	M_{\beta}=\left\lbrace (0,z_2)\in \C^2\mid |z_2|^2=1\right\rbrace
\end{equation*}
which is diffeomorphic to a circle. In particular $M_{\beta}$ is the $S^1$-orbit of the point $(0,1)\in S^3$. The $S^1$-invariant open subset $$U=\left\lbrace (z_1,z_2)\in S^3\mid |z_1|^2<\frac{2}{3}\right\rbrace$$
 is an $S^1$-invariant tubular neighbourhood of the minimal orbit-type stratum, and is diffeomorphic to a solid torus (cf. Figure \ref{hopf}). Since $M_{\beta}$ is an $S^1$-orbit, $U$ is an $S^1$-tubular open subset. We may choose $\mathcal{U'}=\lbrace U\rbrace$. This cover satisfies \ref{minimal ii} of Definition \ref{minimal tubular cover} but it does not satisfy \ref{minimal i}, since it does not cover $M$. To cover $M$ we require the additional open subset
$$V=\left\lbrace (z_1,z_2)\in S^3\mid |z_1|^2>\frac{1}{3}\right\rbrace$$ which is also a solid torus, understood as an $S^1$-invariant tubular neighbourhood of the $S^1$-orbit of $(1,0)$ (cf. Figure \ref{hopf}). It is therefore $S^1$-categorical and then $\cat_{S^1}\left(M\right)\leq 2$. There is in fact equality because otherwise it would mean that $S^3$ is contractible onto a circle, which is untrue.
\vspace{0.5cm}
\begin{figure}[!ht]
	\centering	
		\begin{minipage}[t]{9cm}
		\centering
		\includegraphics[height=6cm]{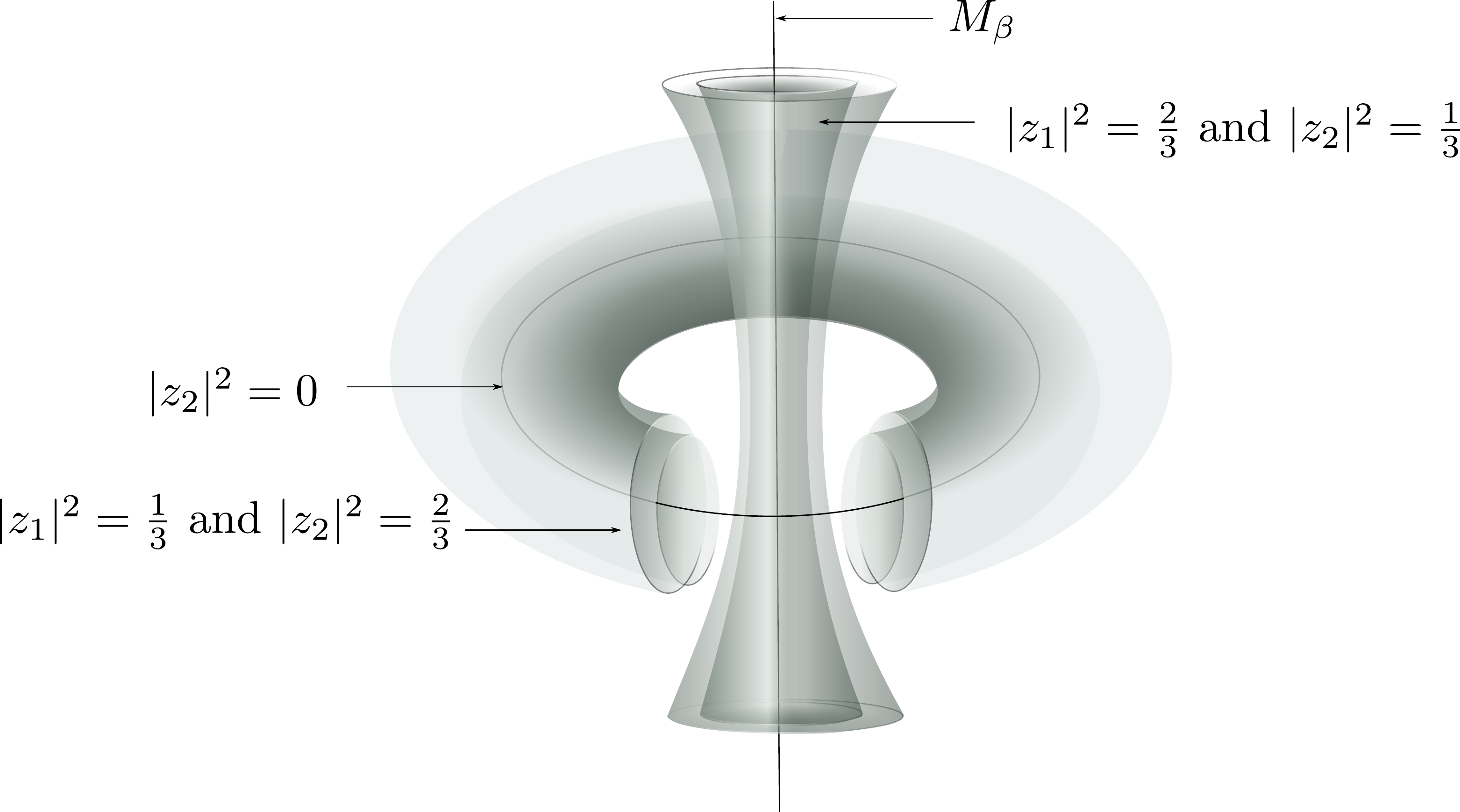}	
		\caption{\small Representation in $\R^3$ of the sphere $S^3$ with a point removed. The stratum $M_{\beta}$ is a circle closing at infinity and the tori around it form a solid torus, which is a tubular neighbourhood.}\label{hopf}
		
\end{minipage}	
\end{figure}
\end{example}

\section{Localization Formula}\label{s: loc formula}

In this section we obtain a localization formula (cf. Corollary \ref{localization thm}) for proper $G$-manifolds which admit a minimal $G$-tubular cover. This formula says in particular that the equivariant LS-category of a proper $G$-manifold is intrinsic to the equivariant LS-category of its minimal orbit-type strata. The theorem below holds in general, without any assumption on the proper $G$-manifold.
\begin{theorem}\label{localization thm first}
		Let $(M,G)$ be a proper $G$-manifold and write $M$ as the disjoint union of its orbit-type strata $\lbrace M_{\beta}\mid \beta\in\mathcal{B}\rbrace$. Let $\mathcal{B}'$ be the biggest subset of $\mathcal{B}$ such that $M_{\beta}$ is minimal for every $\beta\in\mathcal{B}'$. Then $$\cat_G(M)\geq\sum_{\beta\in\mathcal{B}'} \cat_{G}(M_{\beta}).$$
\end{theorem}
\begin{Proof}
	Let $\mathcal{U}$ be a $G$-tubular cover of $M$. Choose $U\in\mathcal{U}$ such that $U\cap M_{\beta}\neq \varnothing$ for some $\beta\in \mathcal{B}'$, say $\beta=((H),(b))$. By Proposition \ref{atmost}, $U$ does not intersect any other minimal stratum and the $G$-deformation retract $H:U\times [0,1]\to U$ retracts onto an orbit $G\cdot x$ of some $x\in M_{\beta}$. The set $V_{\beta}=U\cap M_{\beta}$ is open in $M_{\beta}$ for the subset topology, and it is $G$-invariant because so are $U$ and $M_{\beta}$. 
	
	Let $G\times_{G_x}N_0$ be the local model for $U$. Given $y\in V_{\beta}$ we define the $G$-equivariant path $y(t)=H_t(y)$, where $t\in [0,1]$. In the local model, $y$ reads $[g,\nu]$, and $y(t)$ reads $[g,\nu_t]$ where $\nu_t=(1-t)\nu$. Since $(G_x)_{\nu}\in (H)$, we use the linearity of the $G_x$-action on $N_0$, to obtain $(G_x)_{\nu_t}=(G_x)_{\nu}=(H)$ for all $t\in [0,1]$. Hence $$G_{[g,\nu_t]}=g(G_x)_{\nu_t}g^{-1}\in (H)\quad\mbox{for all}\quad t\in [0,1].$$ In particular, $y(t)\in M_{\beta}$ for all $t\in [0,1]$. Because $y\in V_{\beta}$ is arbitrary and $[0,1]$ is compact, the map $F:V_{\beta}\times [0,1]\to V_{\beta}$ given by $F_t(y)=y(t)$ is a homotopy. It is clearly $G$-equivariant by construction and defines a $G$-deformation retract of $V_{\beta}$ onto $G\cdot x$. It follows that $V_{\beta}$ is $G$-categorical. 
	
	Let $\mathcal{U}_{\beta}\subset\mathcal{U}$ be the subset of all $U\in\mathcal{U}$ such that $U\cap M_{\beta}\neq\varnothing$. Then $$\mathcal{V}_{\beta}=\lbrace V_{\beta}=U\cap M_{\beta}\mid U\in\mathcal{U}_{\beta}\rbrace$$ is a cover of $M_{\beta}$ by $G$-categorical open subsets, which is not necessarily a minimal cover. This procedure associates to each $\beta\in\mathcal{B}'$ a cover $\mathcal{V}_{\beta}$ of $M_{\beta}$.

	Proposition \ref{atmost} says that, if $\alpha,\beta\in\mathcal{B}'$ are distinct, then $\mathcal{U}_{\alpha}\cap\mathcal{U}_{\beta}=\varnothing$. In particular, each $V_{\beta}\in\mathcal{V}_{\beta}$ is determined by a unique $U\in\mathcal{U}_{\beta}$.
	Therefore $$\cat_G(M)\geq \sum_{\beta\in\mathcal{B}'} \cat_{G}(M_{\beta}).$$
	\end{Proof}
Theorem \ref{localization thm first} had already been obtained by\cit{Hurder and T\"oben} (cf. \cite{MR2492297} Theorem $3.7$). However if $(M,G)$ admits a minimal $G$-tubular cover, the following occurs:

\begin{cor}[Localization Formula]\label{localization thm}
		Let $(M,G)$ be a proper $G$-manifold which admits a minimal $G$-tubular cover. Decompose $M$ into its orbit-type strata $\lbrace M_{\beta}\mid\beta\in\mathcal{B}\rbrace$. Let $\mathcal{B}'$ be the biggest subset of $\mathcal{B}$ such that $M_{\beta}$ is minimal for every $\beta\in\mathcal{B}'$. Then $$\cat_G(M)=\sum_{\beta\in\mathcal{B}'} \cat_{G}(M_{\beta}).$$
\end{cor}	
	\begin{Proof}
		By Theorem \ref{localization thm first}, $\cat_G(M)\geq \sum_{\beta\in\mathcal{B}'} \cat_{G}(M_{\beta}).$ The other inequality is a direct consequence of the properties of a minimal $G$-tubular cover (cf. Definition \ref{minimal tubular cover}).
	\end{Proof}
	\begin{proposition}\label{cat norm}
		Let $(M,G)$ be a proper $G$-manifold which admits a minimal $G$-tubular cover. Assume $M_{\beta}$ is a minimal orbit-type stratum with $\beta=((H),(b))$. Then $$\cat_G\left(M_{\beta}\right)=\cat_{N_G(H)}\left(M_{H,b}\right).$$
	\end{proposition}
	\begin{Proof}
		Let $\mathcal{U}'$ be a minimal $G$-tubular cover of $M$ and let $M_{\beta}$ be a minimal orbit-type stratum. By definition of $\mathcal{U}'$, the set $\mathcal{V}_{\beta}=\lbrace V_{\beta}=U\cap M_{\beta}\mid U\in\mathcal{U}'\rbrace$ is the smallest cover by $G$-categorical open subsets of $M_{\beta}$, where the topology of $M_{\beta}$ is the subset topology.

		For every $V_{\beta}\in\mathcal{V}_{\beta}$, let $\widetilde{V}_{\beta}=V_{\beta}\cap M_{H,b}$. Then $\widetilde{V}_{\beta}$ is an $N_G(H)$-invariant open subset of $M_{H,b}$, for the subset topology. Let $H:V_{\beta}\times [0,1]\to V_{\beta}$ be a $G$-deformation retract of $V_{\beta}$ onto some orbit $G\cdot x$ of $x\in M_{\beta}$. Then the $N_G(H)$-homotopy $F:\widetilde{V}_{\beta}\times [0,1]\to \widetilde{V}_{\beta}$ defined by
$$F_t=\restr{H_t}{\widetilde{V}_{\beta}}\quad\mbox{for each}\quad t\in [0,1]$$
is an $N_G(H)$-deformation retract of $\widetilde{V}_{\beta}$ onto the orbit $N_G(H)\cdot x$. Therefore the set $\widetilde{\mathcal{V}}_{\beta}=\lbrace \widetilde{V}_{\beta}=V_{\beta}\cap M_{H,b}\mid V_{\beta}\in\mathcal{V}_{\beta}\rbrace$ is a cover of $M_{H,b}$ made of $N_G(H)$-categorical open subsets. This cover is minimal by assumption and because $M_{\beta}=G\cdot M_{H,b}$. We thus get $$\cat_G\left(M_{\beta}\right)=\cat_{N_G(H)}\left(M_{H,b}\right).$$
	\end{Proof}
	The reader is invited to compare the above result with \cite{marzantowicz} (Proposition $2.1$).
\begin{example}
\normalfont 

\begin{enumerate}[label=(\roman*)]
We verify Theorem \ref{localization thm} on the examples discussed above.
\item Let $M=S^2\subset \R^3$ on which $S^1$ acts by rotations about the $z$-axis. The minimal strata have orbit-type $(H)$ where $H=S^1$, namely
\begin{equation*}
	M_{\beta_b}=\left\lbrace (0,0,(-1)^{b-1})\right\rbrace,
	\quad \beta_b=((S^1),(b))\quad\mbox{and}\quad b=1,2.
\end{equation*}
Then
\begin{equation*}
	\cat_{S^1}\left(M_{\beta_1}\right)+\cat_{S^1}\left(M_{\beta_2}\right)=1+1=\cat_{S^1}\left(M\right).
\end{equation*}
\item Let $M=\C P^2$ equipped with the action of $S^1$
	$$\theta\cdot [z_0:z_1:z_2]=[e^{i\theta}z_0:z_1:z_2].$$ The minimal strata $M_{\beta_1}$ and $M_{\beta_2}$ have orbit-type $(H)=(S^1)$. There is a $\C P^1$ and the single point $[1:0:0]$, respectively. Therefore
$$\cat_{S^1}\left(M_{\beta_1}\right)+\cat_{S^1}\left(M_{\beta_2}\right)=2+1=\cat_{S^1}\left(M\right).$$
 
\item Let $M=S^2$ acted on by the group $\mathbb{T}$ as in Example \ref{ex:TETRAHEDRON} \ref{tetrahedron}. There are three minimal orbit-type strata. Two of them, $M_{\beta_1}$ and $M_{\beta_2}$, have orbit-type $(\Z_3)$. The last minimal stratum $M_{\alpha}$ has orbit-type $(\Z_2)$. We find
 \begin{equation*}
 	\cat_{\Z_3}\left(M_{\beta_1}\right)+\cat_{\Z_3}\left(M_{\beta_2}\right)+\cat_{\Z_2}\left(M_{\alpha}\right)=1+1+1=\cat_{\mathbb{T}}\left(M\right)
 \end{equation*}
	\end{enumerate}
\end{example}

\section{Tubular covers of symplectic toric manifolds}\label{s: tubular toric}

 In this section we show that symplectic toric manifolds admit a minimal tubular cover. Such a cover is constructed explicitly in Theorem \ref{thm toric}. Let $\T$ be torus with Lie algebra $\mathfrak{t}$ and dual Lie algebra $\mathfrak{t}^*$. The smooth action of $\mbox{T}$ on a symplectic manifold $(M,\omega)$ is \defn{Hamiltonian} if there exists a momentum map $\J_{\mbox{\tiny T}}:M\to \mathfrak{t}^*$. The quadruple $(M,\omega,\mbox{T}, \J_{\mbox{\tiny T}})$ is called a \defn{Hamiltonian $\mbox{T}$-manifold}.
 
 \begin{definition}
 	\normalfont Let $\T$ be an $n$-dimensional torus. A Hamiltonian $\T$-manifold $(M,\omega,\T,\J_{\mbox{\tiny T}})$ is called a \defn{symplectic toric manifold} if $(M,\omega)$ is a $2n$-dimensional compact connected symplectic manifold and the Hamiltonian action of $\T$ on $M$ is effective. 
 \end{definition}
For symplectic toric manifolds, the image $\Phi_{\mbox{\tiny T}}(M)$ of the momentum map is a \defn{Delzant polytope} i.e. a convex polytope $\Delta\subset (\R^n)^*$ which is \defn{simple} i.e. each vertex $x$ meets exactly $n$ edges, \defn{rational} i.e. the edges meeting at a vertex $x$ are of the form $x+t\alpha_{x,i}$ where $\alpha_{x,i}\in (\Z^n)^*$, \defn{smooth} i.e. for each vertex $x$ the isotropy weights $\alpha_{x,1},\dots,\alpha_{x,n}$ form a $\Z$-basis of $(\Z^n)^*$.

This observation is due to \cit{Delzant} (cf. \cite{delzant} Lemmas $2.2$ and $2.4$). Delzant also proved that $\Delta$ determines entirely the symplectic toric manifold $(M,\omega,\T,\J_{\mbox{\tiny T}})$, up to $\T$-equivariant symplectomorphisms (cf. \cite{delzant} Theorem $2.1$). His proof relies on a well-known result of convexity obtained independently by\cit{Atiyah} \cite{MR642416} and\cit{Guillemin and Sternberg} \cite{MR664117}, which states that the image of a momentum map for the action of a torus (not necessarily effective) on a compact connected symplectic manifold is a convex polytope.

We recall some standard facts about Morse theory applied to a symplectic toric manifold $(M,\omega,\T,\Phi_{\mbox{\tiny T}})$. The reader is referred to the book of\cit{Guillemin and Sjamaar} (cf. \cite{MR2252111} Section $3.6$) for details. Let $M_{\mbox{\tiny T}}$ be the fixed point set of $\T$. For every $m\in M_{\mbox{\tiny T}}$, the torus acts on the tangent space at $m$. There is a $\T$-invariant complex structure on $M$ such that $T_mM$ is a complex $\T$-representation
with weight space decomposition $$\C_{\alpha_{m,1}}\oplus \dots\oplus \C_{\alpha_{m,n}}$$ where $\alpha_{m,1},\dots,\alpha_{m,n}\in\mathfrak{t}^*$ are the weights of the representation. A \defn{generic component} of the momentum map $\Phi_{\mbox{\tiny T}}:M\to\mathfrak{t}^*$ is a component $\phi^{\xi}=\langle \Phi_{\mbox{\tiny T}}(\cdot),\xi\rangle$ where $\xi\in\mathfrak{g}$ is generic i.e. $\alpha_{m,i}(\xi)\neq 0$ for every $m\in M_{\mbox{\tiny T}}$ and $i=1,\dots,n$. In this case, the critical points of $\phi^{\xi}$ are isolated and $\phi^{\xi}$ is a Morse function whose critical set is precisely $M_{\mbox{\tiny T}}$. Moreover every critical point of $\phi^{\xi}$ has even index. Therefore symplectic toric manifolds possess an extra structure given by the properties of the $\T$-action. This structure is used to construct a minimal $\T$-tubular cover of $M$.

\begin{theorem}\label{thm toric}
Let $(M,\omega,T,\J_{\mbox{\tiny T}})$ be a symplectic toric manifold. Then $M$ admits a minimal $\T$-tubular cover.
\end{theorem}
\begin{Proof}
	Let $\lbrace M_{\beta}\mid \beta\in\mathcal{B}\rbrace$ be the $\mathcal{B}$-stratification of $M$ into orbit-type strata, with strict partial order \eqref{partial order orbit}. Since $\T$ is compact, there are only finitely many minimal orbit-type strata $M_{\beta_1},\dots,M_{\beta_{\ell}}$. By assumption on the $\T$-action, each $M_{\beta_i}$ is an isolated fixed point $m_i\in M_{\mbox{\tiny T}}$. Then there is $\xi_i\in\mathfrak{t}$ such that $\phi^{\xi_i}$ is a generic component of the momentum map, which takes its minimum at $m_i$. Let $-\nabla\phi^{\xi_i}$ be the gradient vector field associated to this component, with corresponding flow $\varphi_t$. Since the image of the momentum map $\Phi_{\mbox{\tiny T}}(M)$ is a Delzant polytope $\Delta$, the $\mathcal{B'}$-stratification $\lbrace \overset{\circ}{F}_{\beta'}\mid \beta'\in\mathcal{B}'\rbrace$ of $\Delta$ by open faces (cf. Section \ref{stratification polytope}) coincides with the $\mathcal{B}$-stratification by orbit-type of $M$. In other words, for every $i=1,\dots,\ell$, we can associate to $\beta_i\in\mathcal{B}$ a unique index $\beta'_i\in\mathcal{B}'$ such that $\Phi_{\mbox{\tiny T}}(M_{\beta_i})$ is precisely the zero-dimensional face $\overset{\circ}{F}_{\beta'_i}$. For each other index $\alpha\in\mathcal{B}$ there is a unique $\alpha'\in\mathcal{B}'$ such that $\Phi_{\mbox{\tiny T}}(M_{\alpha})=\overset{\circ}{F}_{\alpha'}$. Define an open subset $V_{\beta'_i}\subset \mathfrak{t}^*$ by
	\begin{equation*}
		V_{\beta'_i}=\bigcup_{\beta'_i\preceq\alpha'}\overset{\circ}{F}_{\alpha'}.
	\end{equation*}
	By continuity and $\T$-invariance of $\Phi_{\mbox{\tiny T}}$, the subset $U_{\beta_i}=\Phi_{\mbox{\tiny T}}^{-1}(V_{\beta'_i})$ is a $\T$-invariant open neighbourhood of $m_i$ in $M$. It reads
	\begin{equation*}
		U_{\beta_i}=\bigcup_{\beta_i\preceq\alpha}M_{\alpha}.
	\end{equation*}
	For every $m\in U_{\beta_i}\setminus \lbrace m_i\rbrace$, the flow line $\varphi_t(m)$ is defined for every $t\in\R$, by compacity of $M$. By construction of $U_{\beta_i}$, the point $m$ belongs to some orbit-type stratum $M_{\alpha}$ with $\beta_i\prec\alpha$. Since $\varphi_t$ is stratum-preserving, $\varphi_t(m)\in U_{\beta_i}$ for every $t\in\R$. Moreover the only critical point of $\phi^{\xi_i}$ in $U_{\beta_i}$ is $m_i$, and it is a minimum. Hence $\varphi_t(m)$ tends to $m_i$ as $t$ tends to infinity. Therefore the continuous map
\begin{eqnarray*}
		f_m:[0,1[ & \longrightarrow & U_{\beta_i}\\
		t & \longmapsto & \varphi_{\frac{t}{1-t}}(m)
\end{eqnarray*}
extends by continuity into a map $\widetilde{f}_m:[0,1]\to U_{\beta_i}$ with $\widetilde{f}_m(1)=m_i$. Then the map 
\begin{eqnarray*}
		H:U_{\beta_i}\times [0,1] & \longrightarrow & U_{\beta_i}\\
		(m,t) & \longmapsto & \widetilde{f}_m(t)
\end{eqnarray*}
	is a $\T$-deformation retract of $U_{\beta_i}$ onto the orbit $\T\cdot m_i=m_i$. In particular $U_{\beta_i}$ is $\T$-categorical for every $i=1,\dots,\ell$. It is clear that $\mathcal{U}=\lbrace U_{\beta_i}\rbrace_{i=1}^{\ell}$ is a cover of $M$ made of $\T$-tubular open subsets, which are themselves tubular neighbourhoods of the closed strata. By Proposition \ref{atmost}, this cover is the smallest that we can take and then $\mathcal{U}$ is minimal.
\end{Proof}
The result below, due to \cit{Bayeh and Sarkar} (cf. \cite{MR3422738} Theorem $5.1$), is then a direct consequence of the Localization Formula.
\begin{cor}[\cite{MR3422738} Theorem $5.1$]\label{corollary toric}
	Let $(M,\omega,T,\J_{\mbox{\tiny T}})$ be a symplectic toric manifold. Then $\cat_{\mbox{\tiny T}}(M)$ coincides with the cardinality of $M_{\mbox{\tiny T}}$.
\end{cor}
Our choice to consider symplectic toric manifolds makes the proof of Theorem \ref{thm toric} relatively straightforward for two reasons. The first reason is that the fixed points of the $\T$-action are isolated, and the second reason is that the stratification by orbit-type strata of $M$ coincides with the stratification by open faces of the polytope.

\bibliography{LScategory}
\bibliographystyle{plain}

\vspace{0.5cm}
\noindent {\tt School of Mathematics \\
University of Manchester \\
Oxford Road \\
M13 9PL} \\

\vspace{0.5cm}
\noindent {\tt marine.fontaine.math@gmail.com}\\
\noindent {\tt j.montaldi@manchester.ac.uk}
\end{document}